\documentclass[10pt,letter]{amsart}
\usepackage{amssymb,amsbsy,amsmath,amsfonts,amssymb}
\usepackage{latexsym,euscript,exscale}

\title[The complexity of classifying separable Banach spaces]{The complexity of classifying separable Banach spaces up to isomorphism}
\author {Valentin Ferenczi, Alain Louveau, and Christian Rosendal}
\date {June 2006}
\newcommand {\ca} {{2^\N}}
\newcommand{\cantor}{2^{<\omega}}
\newcommand{\baire}{\om^{<\omega}}

\newcommand {\N}{\mathbb N}
\newcommand {\Q}{\mathbb Q}
\newcommand {\R}{\mathbb R}

\renewcommand{\leq}{\ensuremath{\leqslant}}
\renewcommand{\geq}{\ensuremath{\geqslant}}

\newcommand{\tno}{|\!|\!|}
\newcommand{\norm}[1]{\|{#1}\|}

\newcommand{\om}{\omega}

\newcommand{\iso}{\cong}

\newcommand{\tom} {\emptyset}

\newcommand{\hviss}{\leftrightarrow}

\newcommand{\saa}{\Rightarrow}

\newcommand{\til}{\rightarrow}

\newcommand {\Del}{ \; \Big| \;}
\newcommand {\del}{ \; \big| \;}

\newcommand {\go} {\mathfrak}
\newcommand {\ku} {\mathcal}

\newcommand {\e} {\exists}
\renewcommand {\a} {\forall}

\newcommand{\fed}{\boldsymbol}

\newtheorem{thm}{Theorem}
\newtheorem{cor}[thm]{Corollary}
\newtheorem{lemme}[thm]{Lemma}
\newtheorem{prop} [thm] {Proposition}
\newtheorem{defi} [thm] {Definition}

\newtheorem{prob}[thm]{Problem}
\newtheorem{quest}[thm]{Question}

\begin{document}

\thanks{The third author is partially
supported by NSF grant DMS 0556368}


\maketitle

\dedicatory{To Alekos Kechris, on the occasion of his 60th
birthday\footnote{The main result of this paper was obtained at the
end of a conference held in Urbana-Champaign in honour of A.S.
Kechris. We take this opportunity to thank the organisers C. W.
Henson and S. Solecki for their hospitality.}}

\begin{abstract}It is proved that the relation of isomorphism between separable Banach
spaces is a complete analytic equivalence relation, i.e., that any
analytic equivalence relation Borel reduces to it. Thus, separable
Banach spaces up to isomorphism provide complete invariants for a
great number of mathematical structures up to their corresponding
notion of isomorphism. The same is shown to hold for (1) complete
separable metric spaces up to uniform homeomorphism, (2) separable
Banach spaces up to Lipschitz isomorphism, and (3) up to
(complemented) biembeddability, (4) Polish groups up to topological
isomorphism, and (5) Schauder bases up to permutative equivalence.
Some of the constructions rely on methods recently developed by S.
Argyros and P. Dodos.
\end{abstract}

\tableofcontents

\section{Introduction}

A general mathematical problem is that of classifying one class of
mathematical objects by another, that is, given some class $\ku A$,
e.g., countable groups, and a corresponding notion of isomorphism
one tries to find complete invariants for the objects in $\ku A$ up
to isomorphism. In other words, one tries to assign to each object
in $\ku A$ some other object such that two objects in $\ku A$ have
the same assignment if and only if they are isomorphic. This way of
stating it is however slightly misleading, as, in general, one
cannot do better than assigning isomorphism classes in some other
category $\ku B$, or more precisely, on can make an assignment from
$\ku A$ to $\ku B$ such that two objects in $\ku A$ are isomorphic
if and only if their assignments in $\ku B$ are isomorphic. In this
case, we say that we have classified the objects in $\ku A$ by the
objects of $\ku B$ up to isomorphism. However, in order for this not
to be completely trivial, one would like the assignment itself to be
somehow calculable or explicit. The classification should not just
rely on some map provided by the axiom of choice.

For a period going back at least 20 years, there has been a
concentrated effort in descriptive set theory to make a coherent
theory out of the notion of {\em classification} and to determine
which classes of objects can properly be said to be classifiable by
others. The way in which this has been done is by considering
standard Borel spaces that can be considered to fully represent the
classes of objects in question, and then study the corresponding
notion of isomorphism as an equivalence relation on the space. As
standard Borel spaces are fully classified by their cardinality,
which can be either countable or $2^{\aleph_0}$, the perspective
changes from the objects in question to the equivalence relation
instead. One therefore talks of classifying equivalence relations by
each other instead of the corresponding objects. If an equivalence
relation is classifiable by another, one says that the former is
less {\em complex} than the latter. Here is the
precise definition.
\begin{defi}
Let $E$ and $F$ be equivalence relations on standard Borel spaces
$X$ and $Y$ respectively.  We say that $E$ is {\em Borel reducible}
to $F$ if there is a Borel function $f:X\til Y$ such that
$$
xEy\hviss f(x)Ff(y)
$$
for all $x,y\in X$. We denote this by $E\leq_B F$ and informally say
that $E$ less complex than $F$. If both $E\leq_BF$ and $F\leq_BE$,
then $E$ and $F$ are called {\em Borel bireducible}, written
$E\sim_BF$.
\end{defi}
If one looks at the classes of objects that readily are considered
as a standard Borel space $X$, e.g., countable combinatorial and
algebraic objects, separable complete metric structures, one notices
that the corresponding notion of isomorphism is most often analytic,
if not Borel, seen as a subset of $X^2$. Because of this, and also
because the structure theory of $\leq_B$ breaks down beyond the
level of analytic or Borel, one has mostly only developed the theory
in this context.

Classical examples of classifications that fit nicely into this
theory are the classification of countable boolean algebras by
compact metric spaces up to homeomorphism by Stone duality and the
Orstein classification of Bernoulli automorphisms by entropy.

An easy fact, first noticed by Leo Harrington, is that among the
analytic equivalence relations there is necessarily a maximum one
with respect to the ordering $\leq_B$. However, for a long period no
concrete example of this maximum one was found, only abstract
set theoretical versions were known. This problem was solved by Louveau and
Rosendal in \cite{lr}, but at a certain expense. It was noticed that
the definition of the Borel reducibility ordering extends ad
verbatim to quasiorders, i.e., transitive and reflexive relations,
and that one again has a maximum analytic quasiorder.  Using a
representation result that we shall come back to later, it was shown
that, e.g., the relation of embeddability between countable graphs
is a {\em complete} analytic quasiorder, i.e., that its
$\leq_B$-degree is maximum among analytic quasiorders. A simple
argument then shows that the corresponding equivalence relation of
{\em bi}-embeddability is a complete analytic equivalence relation.

A large theory has now been developed concerning analytic and Borel
equivalence relations, but, of course, a main interest in this
theory comes from the fact that it should provide an understanding
of concrete mathematical examples. Thus, in functional analysis,
much effort has been made on trying to understand the structure of
Banach spaces by making inroads into the classification problem,
i.e., by trying to classify separable Banach spaces up to (linear)
isomorphism. Since this is obviously an immensely complicated task
(exactly how immense should be clear from the main result of this
paper), one hoped for a long time that one should instead be
able to find simple subspaces present in every space. However,
even this has turned out somewhat harder than hoped for by several
bad examples of spaces due to Tsirelson \cite{tsi} and, especially,
Gowers and Maurey \cite{gm}.

We will here show exactly how complicated the task is by showing
that the relation of isomorphism between separable Banach
spaces is actually complete as an analytic equivalence relation and
therefore that the classification problem for separable Banach
spaces is at least as complicated as almost any other classification
problem of analysis.

\section{Notation and concepts of descriptive set theory}
In all of the following we will write $\om$ for the set of natural
numbers $\N=\{0,1,2,\ldots\}$ and $2$ for the two-element set
$\{0,1\}$. If $A$ is a non-empty set, a {\em tree} $T$ on $A$ will
be a set of finite strings $t=(a_0,a_1,\ldots, a_{n-1})\in A^{<\om}$
of elements of $A$, $n\geq 0$, containing the empty string $\tom$
and such that if $s\subseteq t$ and $t\in T$, then $s\in T$. Here
$s\subseteq t$ denotes that $s$ is an initial segment of $t$.

A {\em Polish space} is a separable topological space whose topology
is given by a complete metric. The {\em Borel} sets in a Polish
space are those sets that belong to the smallest $\sigma$-algebra
containing the open sets. A {\em standard Borel space} is the
underlying set of a Polish space equipped with the Borel algebra. By
a theorem of Kuratowski, all uncountable standard Borel spaces are
Borel isomorphic with $\R$. An {\em analytic} or $\fed
\Sigma^1_1$-set is a subset of a standard Borel space that is the
image by a Borel function of another standard Borel space. A set is
{\em coanalytic} if its complement is analytic.

A very useful way of thinking of Borel and analytic sets, which is
now known as the Kuratowski--Tarski algorithm, is in terms of the
quantifier complexity of their definitions. Thus, Borel sets are
those that can be inductively defined by using only countable
quantifiers, i.e., quantifiers over countable sets, while analytic
sets are those that can be defined using countable quantifiers and a
single positive instance of a quantifier over a Polish space.

\section{The standard Borel space of separable Banach spaces}
In order to consider the class of separable Banach spaces as a
standard Borel space, we take a separable metrically universal
Banach space, e.g., $X=C([0,1])$ and denote by $F(X)$ the set of all
its closed subsets. We  equip $F(X)$ with its so called
Effros--Borel structure, which is the $\sigma$-algebra generated by
the sets on the form
$$
\{F\in F(X)\del  F\cap U\neq
\tom\},
$$
where $U$ varies over open subsets of $X$. Equipped with this
$\sigma$-algebra, $F(X)$ becomes a standard Borel space, i.e.,
isomorphic as a measure space with $\R$ given its standard Borel
algebra. It is then a standard fact, which is not hard to verify,
that the subset $\go B\subseteq F(X)$ consisting of all the closed
linear subspaces of $X$ is a Borel set in the Effros--Borel
structure. So, in particular, $\go B$ is itself a standard Borel
space, and it therefore makes sense to talk of Borel and analytic
classes of separable Banach spaces, referring by this to the
corresponding subset of $\go B$. We therefore consider $\go B$ as
{\em the} space of separable Banach spaces. It is an empirical fact
that any other way of defining this leads to equivalent results.

The same construction can be done for complete separable metric
spaces. In this case, we begin with a separable complete metric
space, universal for all complete separable metric spaces, for
concreteness we take the Urysohn metric space $\mathbb U$, and then
let $\go M$ be the standard Borel space of all its closed subsets
equipped with the Effros--Borel structure. Again, we see $\go M$ as
the space of all complete separable metric spaces. Since there is a
natural inclusion $\go B\subseteq \go M$, it is reassuring to know
that $\go B$ is a Borel subset of $\go M$.

We shall consider several notions of comparison between Banach
spaces that all will turn out to provide analytic relations on $\go
M$.
\begin{defi}
Let $(X,d_X)$ and $(Y,d_Y)$ be metric spaces. We say that
\begin{itemize}
  \item $X$ and $Y$ are {\em Lipschitz isomorphic} if there is a
  bijection $f:X\til Y$ such that for some $K\geq 1$, we have for all $x$, $y$
  $$
  \frac 1Kd_X(x,y)\leq d_Y(f(x),f(y))\leq Kd_X(x,y).
  $$
  \item $X$ is {\em Lipschitz embeddable} into $Y$ if $X$ is
  Lipschitz isomorphic with a subset of $Y$.
  \item $X$ and $Y$ are {\em uniformly homeomorphic} if there is a
  bijection $f:X\til Y$ such that both $f$ and $f^{-1}$ are
  uniformly continuous.
  \end{itemize}
\end{defi}
To see that, for example, the relation of Lipschitz isomorphism on
$\go M$ is analytic, we notice that for $X,Y\in \go M$, $X$ is
Lipschitz isomorphic to $Y$ if and only if
\begin{displaymath}\begin{split}
\e K\geq 1\;\e (x_n)\; \e (y_n)\;\Big( &\a m\;\big(X\cap
U_m\neq\tom\til \e n\; x_n\in U_m\big)\\
\&\; &\a m\;\big( Y\cap U_m\neq \tom\til\e n\; y_n\in
U_m\big)\\
\&\; &\big (\a n\; x_n\in X\big)\; \&\;\big(\a n\; y_n\in
Y\big)\\
\&\; &\a n,m\;\big(\frac 1Kd_X(x_n,x_m)\leq d_Y(y_n,y_m)\leq
Kd_X(x_n,x_m)\big)\Big)
\end{split}\end{displaymath}
where $\{U_m\}$ is an open basis for $\mathbb U$. In other words,
$X$ and $Y$  are Lipschitz isomorphic if and only if they have
Lipschitz isomorphic countable dense subsets and this condition can
be expressed in an analytic manner in $\go M^2$.

\section{Results}
We are now ready to explain the results of the present paper.
The main idea is to combine a refinement of the completeness results
of Louveau and Rosendal in \cite{lr} with a recent construction by
S. Argyros and P. Dodos from \cite{ad} that was done for different
though related purposes. The thrust of the Agyros--Dodos
construction is to be able to associate to each analytic set of
Schauder bases (i.e., analytic set of subsequences of the universal Pe\l
czy\'nski basis) a separable space that essentially only has basic
sequences from the analytic set. Of course this is not quite
possible (take, for example, $\{\ell_1,\ell_2\}$), but at least one
can get a significant amount of control over the types of basic sequences
present. Now the completeness results of \cite{lr}, on the other
hand, show exactly that certain relations related to identity (in
the codes) of analytic sets are complete analytic equivalence
relations, and we are therefore able to code these relations into
relations between Banach spaces.

Our first result concerns the relation of permutative equivalence
between Schauder bases. First, as mentioned above, we take as our
standard Borel space of Schauder bases, the set $[\om]^\om$ of all
infinite subsets of $\om=\{0,1,\ldots\}$, where we identify a subset
of $\om$ with the corresponding subsequence of the universal Schauder basis
constructed by Pe\l zcy\'nski \cite{p}. It was proved in
\cite{r:cof} that the relation of equivalence between Schauder bases
was a complete $\fed K_\sigma$-quasiorder, and here we prove
\begin{thm}\label{perm}
The relation of permutative equivalence between (even unconditional)
Schauder bases is a complete analytic equivalence relation.
\end{thm}
This result gives probably the a priori simplest naturally occurring
equivalence relation which is known to be analytic complete. Thus,
somewhat surprisingly, the relation of isomorphism between separable
Banach spaces is Borel reducible to permutative equivalence. One can
actually consider this even as a sort of representation result for,
e.g., separable Banach spaces. We can in a Borel manner associate to
each space a basis such that two spaces are isomorphic if and only
if the two bases are permutatively equivalent. If this could be done in a more
informative or explicit manner than in our construction, one could
really hope for an increased understanding of the isomorphism
relation in terms of the more readily understandable relation of
permutative equivalence, and, in fact, one could consider such a
result a best positive solution to the problem of representing separable
Banach spaces by bases.

We subsequently use Theorem \ref{perm} to study uniform
homeomorphism between complete separable metric spaces. Though we
have not been able to show that this relation restricted to $\go B$
is complete analytic, we do get
\begin{thm}
The relation of uniform homeomorphism between complete separable
metric spaces is a complete analytic equivalence relation.
\end{thm}
There is of course another perhaps more immediate relation to study
on $\go M$, namely isometry. The situation for this relation is
nevertheless slightly different, for Gao and Kechris \cite{gk} have
shown that this relation is bireducible with the most complex orbit
equivalence relation, $E_G$, induced by the continuous action of a
Polish group on a Polish space. And by results of Kechris and
Louveau \cite{kl}, this relation is strictly less complex than a
complete analytic equivalence relation. So this also means that
isometry on $\go B$ is simpler than permutative equivalence of
bases. Recently, Melleray \cite{mel} has been able to show that also
restricted to $\go B$ is the (linear) isometry relation Borel
bireducible with $E_G$. This should be contrasted with the result in
\cite{lr} saying that the relation of (linear) isometric
biembeddability on $\go B$ is a complete analytic equivalence
relation.

The most important relation between Banach spaces is however the
relation of linear isomorphism, which has turned out to be
exceedingly difficult to understand, so much that among Banach space
theorists there has even been a feeling that the category of Banach
spaces might not be the right category to study, but that one
instead should consider only spaces with a basis. Theorem \ref{perm}
of course shows that such a restriction would not really decrease
the complexity of the task, but at least the following result should
be a comfort in the sense that it confirms the feeling of outmost
complexity.
\begin{thm}
The relations of isomorphism and Lipschitz isomorphism between
separable Banach spaces are complete analytic equivalence
relations.
\end{thm}
This result is the culmination of a series of successive lower
estimates of the complexity by Bossard \cite{b}, Rosendal
\cite{r:these}, and Ferenczi--Galego \cite{fg}. We also consider the
corresponding quasiorders of embeddability etc. and show that these
are also complete in their category.

We finally consider Banach spaces as abelian groups and notice that
any continuous group isomorphism is also linear. Therefore,
\begin{thm}
The relation of topological isomorphism between (abelian) Polish
groups is a complete analytic equivalence relation.
\end{thm}
Previous work on the complexity of isomorphism between groups has
been exclusively on the countable discrete case. An early result of
Friedman and Stanley \cite{fs} states that the relation of
isomorphism between countable discrete groups is complete among all
isomorphism relations between countable structures, while Thomas and
Velickovic \cite{tv} prove that restricted to the class of finitely
generated groups it becomes complete among all Borel equivalence
relations having countable classes.

\section{A variant of the completeness method}
In \cite{lr}, Louveau and Rosendal established a representation
result for  analytic quasiorders, and used this result  to prove that some
$\fed\Sigma^1_1$ quasiorders are {\em complete}, i.e., have the
property that any other $\fed\Sigma^1_1$ quasiorder is Borel
reducible to them, and to deduce from this that certain $\fed\Sigma^1_1$
equivalence relations are also complete .

It was clear,  from the way the completeness results
were derived from the representation  that the technique was flexible and
could lead to improved results.
This was implicitly acknowledged in \cite{lr} and stated more
explicitly in \cite{r:cof} in the case of $\fed
K_\sigma$  quasiorders, but, as at the time no essential use
of it was made, the details were not spelled out. However,
in the applications in the present paper, the finer versions
have turned out to be crucial for our proofs and we
therefore proceed to state the results precisely.

The main idea is to desymmetrise the situation, both for the relations
and for the reducibility ordering.
We think of a binary relation $R$ on some $X$ as the pair of
relations $(R, \neg R)$, where $\neg R$ denotes the complement of $R$ in $X^2$. With this
identification, Borel reducibility is now defined on certain kinds of pairs, and we extend it
to arbitrary pairs as follows.

\begin{defi}
Let $(R_1,R_2)$ and $(S_1,S_2)$ be two pairs of binary relations
on  standard Borel spaces $X$ and $Y$ respectively.
A Borel map $f: X \to Y$ is a {\em Borel homomorphism} from $(R_1,R_2)$
to $(S_1,S_2)$ if for all $x,y \in X$, $xR_1y \rightarrow f(x)S_1f(y)$,
and $xR_2y \rightarrow f(x)S_2f(y)$.
We  say that $(R_1,R_2)$ is {\em Borel hom-reducible} to $(S_1,S_2)$, and
write
$$
(R_1,R_2)\preccurlyeq_B(S_1,S_2)
$$
if there is a Borel homomorphism from $(R_1,R_2)$ to $(S_1,S_2)$.
\end{defi}
Borel hom-reducibility is clearly a quasi-ordering, as homomorphisms
can be composed. Moreover, one has from the definitions that
$$
R\leq_B S \leftrightarrow (R,\neg R) \preccurlyeq_B(S,\neg S),
$$
so that using the identification above, $\preccurlyeq_B$ is
indeed an extension of $\leq_B$.

In the sequel, we will let $(R_1,R_2)\leq_B R$ and $R\leq_B
(R_1,R_2)$ stand for $(R_1,R_2)\preccurlyeq_B (R,\neg R)$ and
$(R,\neg R) \preccurlyeq_B (R_1,R_2)$, respectively.

Suppose now that we are interested in a class $\ku C$ of binary relations on
standard Borel spaces, e.g., analytic quasiorders or analytic equivalence
relations. We  say that a pair $(R_1,R_2)$ of binary relations
is $\ku C$-{\em hard} if any element of $R\in\ku C$ is Borel hom-reducible to
$(R_1,R_2)$, or more precisely $R\leq _B(R_1,R_2)$. And we say that $(R_1,R_2)$ is $\ku C$-complete if it is $\ku C$-hard,
and, moreover, it is Borel hom-reducible to some element of $\ku C$.
In other words, if we set $$\ku C^* = \{(R_1,R_2)\del  \exists R \in \ku  C ~~(R_1,R_2)\leq_B R\},$$
the pair $(R_1,R_2)$ is $\ku C$-complete if it is $\preccurlyeq_B$-maximum in $\ku C^*$.
It is very easy to check that there is a $\ku C$-complete $R$, i.e.,  a $\leq_B$-maximum
element in $\ku C$, if and only if there is a  $\ku C$-complete pair $(R_1,R_2)$, and, moreover,
if this happens, then the $\ku C$-complete $R$'s are exactly the ones which, viewed as
pairs, are complete.

However, the strong completeness results we will need in the sequel only rely  on the
following  simple observation: Suppose $(R_1,R_2)$ is $\ku C$-hard and
$(R_1,R_2)$ reduces to some $R\in \ku C$. Then $R$ is in fact $\ku
C$-complete. This is
of course an obvious fact, but will be quite handy as it will allow
us to work in some cases with a more manageable pair $(R_1,R_2)$
than with a single complete relation $R$.

To give the flavour of the arguments, consider the case of
$\fed\Sigma^1_1$ equivalence relations. Then there is a pair which
is easily seen to be hard for this class: suppose we are given a
coding $\alpha \mapsto A_\alpha$ of $\fed\Sigma^1_1$ subsets of say
$2^\omega$ by elements of some Polish space $X$ (we will be more
specific later on), and define binary relations
$=_{\fed\Sigma^1_1}$,  $\subseteq_{\fed\Sigma^1_1}$, ${\rm
Disj}_{\fed\Sigma^1_1}$ on $X$ corresponding, ``in the codes'', to
$=$, $\subseteq$, and disjointness between non-empty
$\fed\Sigma^1_1$ sets, i.e., set
\[\begin{split}
\alpha =_{\fed\Sigma^1_1}\beta &\leftrightarrow A_\alpha \not=\emptyset \;\&\; A_\beta \not=\emptyset \;\&\; A_\alpha = A_\beta,\\
\alpha \subseteq_{\fed\Sigma^1_1}\beta &\leftrightarrow A_\alpha\not=\emptyset \;\&\; A_\beta \not=\emptyset \;\&\; A_\alpha \subseteq A_\beta,\\
\alpha \;{\rm Disj}_{\fed\Sigma^1_1}\beta &\leftrightarrow
A_\alpha\not=\emptyset \;\&\; A_\beta \not= \emptyset \;\&\;
A_\alpha\cap A_\beta = \emptyset.
\end{split}\]
For any reasonable coding, the pair $(=_{\fed\Sigma^1_1}, {\rm
Disj}_{\fed\Sigma^1_1})$ is hard for the class of $\fed\Sigma^1_1$
equivalence relations. For if $E$ is such a relation, which without
loss of generality we can view as defined on $2^\omega$, one can
associate to $x\in 2^\omega$ a code for its equivalence class
$[x]_E$ in a continuous way, and this gives a homomorphism from $E$
to the pair $(=_{\fed\Sigma^1_1}, {\rm Disj}_{\fed\Sigma^1_1})$. We
do not know if this pair is complete for analytic equivalence
relations, hence if it can be used to obtain completeness results.
But the basic representation of \cite{lr} allows to replace it by a
complete pair.

Let us first recall this result:
\begin{thm}[A. Louveau and C. Rosendal \cite{lr}]\label{representation}
Let $R\subseteq 2^\om\times 2^\om$ be a $\fed\Sigma_1^1$ quasiorder.
Then there exists a tree $T$ on $2\times 2\times \om$ with the
following properties
\begin{enumerate}
  \item\label{1} $xRy\hviss \e \alpha\in \om^\om\; \a n\;
  (x|_n,y|_n,\alpha|_n)\in T$,
  \item\label{2} if $(u,v,s)\in T$ and $s\leq t$, then $(u,v,t)\in T$,
  \item\label{3} for all $(u,s)\in (2\times \om) ^{<\om}$, we have $(u,u,s)\in
  T$,
  \item\label{4} if $(u,v,s)\in T$ and $(v,w,t)\in T$, then $(u,w,s+t)\in T$.
\end{enumerate}
\end{thm}
In the statement of the above theorem, if $s$  is a finite sequence and $|s|$ denotes its length, we let for
sequences $s$ and $t$, $s\leq t$ mean that $|s|=|t|$ and for all $i<|s|,  s(i)\leq t(i)$.  And
for $s,t$ of the same length, we let $(s+t)(i)=s(i)+t(i)$.

Let $\go T$ be the class of non-empty {\em normal} trees on $2\times\om$,
i.e., trees $T$ with $(\emptyset,\emptyset) \in T$ and such that whenever
$(u,s)\in T$ and $s\leq t$, then
also $(u,t)\in T$. Viewed as a subset of $2^{(2\times \om)^{<\om}}$, it is closed,
hence a (compact) Polish space.

We view each normal tree $T$ as coding the $\fed\Sigma^1_1$ set
$$
A(T)=\{\alpha\in 2^\om \del \exists \beta\in\om^\om \;\forall n \;(\alpha |_n,\beta |_n)\in T\}.
$$
As is well known, any $\fed\Sigma_1^1$ subset of $2^\om$ is on the form $A(T)$ for some $T$, so we really have a coding.

\begin{defi}
We define the following binary relations on $\go T$.
For $S,T\in\go T$ we let
\[\begin{split}
S\leq_{\fed\Sigma_1^1}T &\hviss \e \alpha\in \om^\om\;
\a(u,s)\;\big((u,s)\in S \til(u,s+\alpha|_{|s|})\in T\big),\\
S\equiv_{\fed\Sigma^1_1}T &\hviss S\leq_{\fed\Sigma_1^1}T
\;\&\; T\leq_{\fed\Sigma_1^1}S,\\
S\not\subseteq_{\fed\Sigma_1^1}T &\hviss A(S)\not\subseteq A(T),\\
S \;{\rm Disj}_{\fed\Sigma^1_1} T&\hviss A(S)\neq \emptyset\;\&\; A(T) \neq \emptyset\;\&\;
A(S)\cap A(T) =\emptyset,\\
S\neq_{\fed\Sigma^1_1}T &\hviss A(S)\neq A(T).
\end{split}\]
\end{defi}

\begin{thm}

(i) The pair $(\leq_{\fed\Sigma_1^1},\not\subseteq_{\fed\Sigma_1^1})$ is
complete for the class $\ku C_{\rm qo}$ of analytic quasiorders.

(ii) The pair $(\equiv_{\fed\Sigma^1_1}, {\rm Disj}_{\fed\Sigma^1_1})$ is
complete for the class $\ku C_{\rm eq}$ of analytic equivalence relations,
and hence, a fortiori, the
pair $(\equiv_{\fed\Sigma_1^1}, \neq_{\fed\Sigma^1_1})$ is complete for $\ku C_{\rm eq}$ too.
\end{thm}

\begin{proof}

Note first that $\leq_{\fed\Sigma_1^1}$ is an analytic quasiorder, and $\equiv_{\fed\Sigma^1_1}$
an analytic equivalence relation, with trivially
$$
(\leq_{\fed\Sigma_1^1},\not\subseteq_{\fed\Sigma_1^1}) \leq_B \;\leq_{\fed\Sigma_1^1}
$$
and
$$
(\equiv_{\fed\Sigma^1_1}, {\rm Disj}_{\fed\Sigma^1_1})
\preccurlyeq_B (\equiv_{\fed\Sigma_1^1}, \neq_{\fed\Sigma^1_1})
\leq_B \;\equiv_{\fed\Sigma^1_1},
$$
via the identity map, so that it is enough to prove that the pairs are hard for their respective
classes.

For part (i), let $R$ be a $\fed\Sigma_1^1$ quasiorder on some Polish space $X$. Embedding $X$
into $2^\om$ in a Borel way, we may assume $R$ is defined on $2^\om$. Let then $T$ be the tree
given by theorem 8 and define a continuous map $f:2^\om \to \go T$ by
$$
f(x) = \{(u,s) \in (2\times
\om)^{<\om}\del (u,x|_{|u|},s)\in T\}.
$$
We claim this map works. First, each $f(x)$ is indeed a non empty normal tree by properties (2)
and (3) of $T$. Also, by property (1), $A(f(x)) = \{y\in \ca\: yRx\}$, so that if $\neg xRy$, we get
$x\in A(f(x))$ but $x\not\in A(f(y))$, whence $f(x) \not\subseteq_{\fed\Sigma_1^1} f(y)$.
Conversely suppose $xRy$. Then by property (1) of $T$ there is some $\alpha \in \om^\om$ with
$(x|_n, y|_n,\alpha|_n)\in T$ for all $n$. But then this $\alpha$ witnesses $f(x)
\leq_{\fed\Sigma_1^1} f(y)$, for if $(u,s)\in f(x)$, i.e., if $(u,x|_{|u|},s)\in T$, we get from
property (4) of $T$ that $(u,y|_{|u|},s+\alpha|_{|u|}) \in T$, as $(x|_{|u|},
y|_{|u|},\alpha|_{|u|})\in T$. So $(u,s+\alpha|_{|u|}) \in f(y)$, as desired. This proves (i).

To prove (ii), we again assume that the Polish space is $2^\om$. As an equivalence relation  $E$
is in particular a quasiorder, we can apply part (i), and get a continuous map $f:2^\om \to \go T$
such that $xEy \rightarrow f(x) \leq_{\fed\Sigma_1^1} f(y)$ and $A(f(x))=[x]_E$. But then
trivially $f$ hom-reduces $(E,\neg E)$ to $(\equiv_{\fed\Sigma^1_1}, {\rm
Disj}_{\fed\Sigma^1_1})$, as wanted.
\end{proof}

The previous result is conceptually the simplest one. Unfortunately, we will need in the sequel
a slight improvement of the last statement, obtained by restricting the domains of the relations
to {\em pruned} normal trees, which makes things messier.

Recall that a non empty tree is {\em pruned} if any sequence in it admits a strict extension
that is still in it.

Let
$\go T_{\rm pr}$ be the $G_\delta$ subset of $\go T$ consisting  of the non-empty pruned normal trees,
 and denote by $\leq_{\fed\Sigma_1^1}^{\rm pr}$, $\equiv_{\fed\Sigma^1_1}^{\rm
pr}$, $\not\subseteq_{\fed\Sigma_1^1}^{\rm pr}$ and $\neq_{\fed\Sigma^1_1}^{\rm pr}$ the
restrictions to $\go T_{\rm pr }$ of the corresponding relations.

\begin{thm}
(i) The pair $(\leq_{\fed\Sigma^1_1}^{\rm pr},\not\subseteq_{\fed\Sigma_1^1}^{\rm pr})$  is
complete for the class $\ku C_{\rm qo}$.

(ii) The pair $(\equiv_{\fed\Sigma^1_1}^{\rm
pr},\neq_{\fed\Sigma^1_1}^{\rm pr})$ is complete for the class $\ku C_{\rm eq}$.
\end{thm}

\begin{proof}
(ii) follows from (i) as before, and by the previous theorem, it is enough to prove that
$(\equiv_{\fed\Sigma^1_1},\neq_{\fed\Sigma^1_1}) \preccurlyeq_B
(\equiv_{\fed\Sigma^1_1}^{\rm pr},\neq_{\fed\Sigma^1_1}^{\rm pr})$.

Let $T$ be a non-empty normal tree. We define a tree $T^*$ as follows: for each $(u,s)$ in $T$
of length say $n$, put in $T^*$ all sequences $(u',s')$ of length $\geq 2n$ which satisfy

\noindent (a) $\forall i<n \;u'(2i)=u(i) \;\&\; s'(2i)=s(i)$,

\noindent (b) $\forall i<n \;u'(2i+1)=0$,

\noindent and (c) $\forall i\geq 2n \;u'(i)=1$,

\noindent together with their initial segments. Easily $T^*$ is
still normal, and is now pruned as any sequence in $T^*$ can be
extended using (c). So this defines a continuous map from $\go T$ to
$\go T_{\rm pr}$, and it is enough to check it is the homomorphism
we want. For each $\alpha \in \om^\om$, set $\alpha^*(2i)=\alpha(i)$
and $\alpha^*(2i+1)=0$. Then one checks easily using (a) and (b)
that if $\alpha$ is a witness that $S \leq_{\fed\Sigma^1_1} T$, then
$\alpha^*$ witnesses that $S^* \leq_{\fed\Sigma^1_1}^{\rm pr} T^*$.
Let also $D_1 \subseteq \ca$ be the countable set of eventually $1$
sequences, and for $u\in 2^n$, let $\alpha_u(2i)=u(i)$ and
$\alpha_u(2i+1)=0$ for $i<n$, and $\alpha_u(i)=1$ for $i\geq 2n$.
Then easily from (a), (b), and (c) one gets that
$$
A(T^*)=\{\alpha^*\del
\alpha \in A(T)\}\cup \{\alpha_u\del \exists s~~~(u,s)\in T\}.
$$
From this we get that $A(S)\not \subseteq A(T)$ implies that
$A(S^*)\not \subseteq A(T^*)$ as desired.

Note that one does not have necessarily $A(S^*)\cap A(T^*)=\emptyset$ when  $A(S)\cap
A(T)=\emptyset$. Still we could define $S\;{\rm Disj}_{\fed\Sigma^1_1}^{\rm pr}T$ for  $S,T\in \go T_{\rm
pr}$ by
$$
A(S^*)\setminus D_1 \neq \emptyset\;\&\; A(T^*)\setminus D_1 \neq \emptyset\;\&\;
A(S^*)\cap A(T^*)\subseteq D_1
$$
and get a slight improvement on part (ii), as this last relation
is both smaller and descriptively simpler than $\neq_{\fed\Sigma^1_1}^{\rm pr}$.

\end{proof}

The interesting part in this result is the following: if one wants to prove that a
certain analytic equivalence relation $E$ is complete by providing a
reduction from normal (pruned) trees, one needs to show that if
the two trees code the same analytic set in a strong sense, namely
that there is a uniform $\alpha$ that can translate between the
codes, then the images are $E$-equivalent. But on  the other hand,
for the negative direction one only needs to consider trees that
really code different analytic sets and show that their images are
$E$-inequivalent.
In the applications, we will construct objects from normal trees that ``realise only
the types'' given by the analytic set corresponding to the tree.
For example, in the case of separable Banach spaces, we shall
construct from a normal tree coding an analytic subset of $]1,2[$ a
Banach space whose only $\ell_p$ subspaces are exactly $\ell_2$ plus
those given by the analytic set. Thus if two normal trees are
$\neq_{\fed\Sigma^1_1}$ related, they have different $\ell_p$
subspaces and are hence non-isomorphic. A similar line of thinking in
terms of extreme pairs of quasiorders is also present in Camerlo \cite{c}.

Before we go to Banach spaces, let us illustrate the previous discussion
with a natural example of a
$\ku C_{\rm qo}$-complete pair which could potentially be of use
elsewhere. In analogy with the case of separable Banach spaces,
where the class of $\ell_p$ subspaces of a space will turn out to be
sufficient to separate non-isomorphic spaces, we search for simple
types of objects in a certain category and then associate with each
object its ``spectrum'' consisting of the simple types embeddable
into it.

Let $\go A$ be the class of {\em combinatorial} trees on $\N$, i.e.,
acyclic, connected, symmetric relations on $\N$ and let $\go A_{\rm
f}$ be the subclass of trees of finite valency. For each $T\in \go
A$, we let $\sigma(T)$ be the {\em spectrum} of $T$, which is the
set of all $S\in \go A_{\rm f}$ that embed into $T$. We then let
$\sqsubseteq$ be the relation of embeddability between combinatorial
trees and put $S\subseteq_\sigma T$ if $\sigma(S)\subseteq
\sigma(T)$. Using a simple modification of the construction in
\cite{lr} showing that $\sqsubseteq$ is a complete analytic
quasiorder one can show that
\begin{prop}
    The pair $(\sqsubseteq,\not\subseteq_\sigma)$ is $\ku C_{\rm
    qo}$-complete.
\end{prop}

In view of this proposition, it would be interesting to find
combinatorial realisations of the complete analytic equivalence
relation. However, as is well-known, no Borel class of countable
model-theoretical structures provides an isomorphism relation which
is complex enough. Instead, one should look for other notions of
isomorphism. A potential candidate would be the relation of
quasi-isometry between countable graphs. Here two graphs $R$ and $S$
on the vertex set $\N$ are said to be quasi-isometric if there is a
function $\phi:\N\til \N$ and numbers $K$, $N$ and $L$ such that for
all $n,m$
$$
d_R(n,m)\leq Kd_S(\phi(n),\phi(m))+N,
$$
$$
d_S(\phi(n),\phi(m))\leq Kd_R(n,m)+N,
$$
and
$$
\a k\;\e l\;d_S(k,\phi(l))\leq L.
$$
S. Thomas \cite{t} has recently proved that this relation, when
restricted to locally finite graphs, is Borel bireducible with the
complete $\fed K_\sigma$-equivalence relation.

\

To end up this section, let us discuss another situation where it is possible to get a nice
complete pair (although we have no application for it). It is the case of orbit equivalence
relations for Borel actions of Polish groups.

Fix a Polish group $G$. If $X$ is a standard Borel space and $\alpha: G\times X\to X$ is a
Borel action of $G$ on $X$, one defines the associated $\fed\Sigma^1_1$ orbit equivalence
$E_G^X$ by
$$xE_G^Xy \hviss \exists g\in G~~~\alpha(g,x)=y.$$ We let $\ku C_G$ be the class of all such
orbit equivalence relations. By a result of Becker--Kechris \cite{becker}, it is the same class, up to
Borel isomorphism, as the class of orbit equivalences corresponding to continuous actions of
$G$ on Polish spaces $X$.

Becker and Kechris also proved that there is a complete element in $\ku C_G$. We now provide a
complete pair for it (which gives a somewhat different complete element).

First, fix some universal Polish space $X_0$, like the Urysohn space or $\R^\om$, with the
property that any Polish space is homeomorphic to a closed subspace of it. Let $Z$ be the
standard Borel space of non-empty closed subsets of $X_0\times G$, equipped with the Effros--Borel structure.
Define an action $(g,F) \mapsto g.F$ of $G$ on $Z$ by setting
$g.F=\{(x,gh)\del (x,h)\in F\}$. It is easy to check that this action is Borel, and hence the
associated $E_G^Z$ is in $\ku C_G$.

For $F$ in $Z$, set $A(F) = \{x\in X_0\del \exists g\in G~~~(x,g)\in F\}$, and define ${\rm
Disj}_G$ on $Z$  by $F \;{\rm Disj}_G F' \hviss A(F)\cap A(F')=\emptyset$.

\begin{thm}
The pair $(E_G^Z, {\rm Disj}_G)$ is complete for the class $\ku C_G$ (and hence
$E_G^Z$ is complete too).
\end{thm}
\begin{proof}
As $E_G^Z$ is in $\ku C_G$, we only have to check that $(E_G^Z,{\rm Disj}_G)$ is $\ku
C_G$-hard. And by the result of Becker and Kechris quoted above, we only have to consider
Polish spaces $X$ and continuous actions $\alpha: G\times X \to X$. View $X$ as a closed subset
of $X_0$ and associate to each  $x\in X$ the element $F_x \in Z$ defined by $$F_x=\{(y,g)\in
X_0\times G\del y\in X\;\&\; \alpha (y,g)=x\}.$$ Note that $F_x$ is non-empty as $(x,1_G)\in
F_x$, and one can check, using the continuity of $\alpha$, that the map $x\mapsto F_x$ is
Borel. One easily also checks that $F_{\alpha(g,x)}=g.F_x$, so that if $xE_G^Xy$, then
$F_xE_G^ZF_y$. And finally $A(F_x)$ is just the orbit of $x$ for the action $\alpha$, hence if
$\neg xE_G^Xy$, $A(F_x)\cap A( F_y) = \emptyset$, as desired.
\end{proof}

\section{An $\ell_p$-tree basis}
We define in this section the construction of a basic sequence from
a tree on $2\times \om$. This will prove to be fundamental in our
later proofs. We begin by choosing a Cantor set of $p$'s in the
interval $]1,2[$ in the following fashion. The set is given by a
Cantor scheme $(I_u)_{u\in \cantor}$ of non-empty closed
subintervals $I_u\subseteq ]1,2[$ such that the following holds.
\begin{enumerate}
\item $I_{u0}\cup I_{u1}\subseteq I_u$,
\item $\max I_{u0}<\min I_{u1}$,
\item $I_{u0}$ contains the left endpoint of $I_u$,
\item $I_{u1}$ contains the right endpoint of $I_{u}$,
\item\label{eq-int} the standard unit vector bases of $\ell_{\min I_u}^{|u|}$ and $\ell_{\max
I_u}^{|u|}$ are $2$-equivalent, i.e., $\frac{|u|^{1/{\min I_u}}}{|u|^{1/{\max
I_u}}}\leq2$, for all $u\neq \tom$.
\end{enumerate}
If now $\alpha\in \ca$, we denote by $p_\alpha$ the unique point in
$\bigcap_{u\subseteq \alpha}I_u$. Then
$$
\alpha\in \ca\mapsto
p_\alpha\in ]1,2[
$$
is an orderpreserving homeomorphism between $\ca$ with the
lexicographical ordering and a compact subset of $]1,2[$.

In the following we denote by $\mathsf T$ the complete normal tree
$(2\times \om)^{<\om}$. As always we identify the elements of
$\mathsf T$ with the pairs $t=(u,s)\in \cantor\times \baire$ such
that $|u|=|s|$. A {\em segment} $\go s$ of $\mathsf T$ is just a set
on the form $\go s=\{t\in\mathsf T\del t_0\subseteq t\subseteq
t_1\}$ for some $t_0, t_1\in \mathsf T$.

We now let $\ku V=c_{00}(\mathsf T)$ be the vector space with basis
$(e_t)_{t\in \mathsf T}$. For each segment
$$\go s=\{(u_0,s_0)\subsetneq (u_1,s_1)\subsetneq
\ldots\subsetneq (u_n,s_n)\}
$$
of $\mathsf T$, we define a norm $\|\cdot \|_\go s$ on $\ku V$ as
follows
$$
\|\sum_{t\in \mathsf T}\lambda_te_t\|_\go s=\sup_{m\leq
n}\|(\lambda_{(u_0,s_0)},\lambda_{(u_1,s_1)},\ldots,\lambda_{(u_m,s_m)})\|_{\min
I_{u_m}}.
$$
We notice that for $m\leq n$, we have $u_m\subseteq u_n$ and so
$I_{u_n}\subseteq I_{u_m}$, whence by condition (\ref{eq-int})
\begin{displaymath}
\begin{split}
&\|(\lambda_{(u_0,s_0)},\lambda_{(u_1,s_1)},\ldots,\lambda_{(u_m,s_m)})\|_{\min
I_{u_m}} \\
\leq&
2\;\|(\lambda_{(u_0,s_0)},\lambda_{(u_1,s_1)},\ldots,\lambda_{(u_m,s_m)})\|_{\min
I_{u_n}} \\
 \leq&
2\;\|(\lambda_{(u_0,s_0)},\lambda_{(u_1,s_1)},\ldots,\lambda_{(u_m,s_m)},\ldots,\lambda_{(u_n,s_n)})\|_{\min
I_{u_n}}.
\end{split}
\end{displaymath}
Thus, if $\sigma=(\alpha,\beta)\in [\mathsf T]$ is a branch of
$\mathsf T$ containing the segment $\go s$, then $p_\alpha\in
I_{u_n}$, and thus
\begin{equation}\label{eq with p}
\|\sum_{t\in\go s}\lambda_te_t\|_{p_\alpha}\leq\|\sum_{t\in\go
s}\lambda_te_t\|_{\go s} =
\|\sum_{t\subseteq\sigma}\lambda_te_t\|_\go s\leq 2 \|\sum_{t\in\go
s}\lambda_te_t\|_{p_\alpha}\leq 2 \|\sum_{t\subseteq
\sigma}\lambda_te_t\|_{p_\alpha}
\end{equation}
Finally, we define the norm $\tno \cdot\tno $ on $\ku V$ by
$$
{\tno }\sum_{t\in \mathsf T}\lambda_te_t{\tno }=\sup \big\{\big
(\sum_{i=1}^l\big\|\sum_{t\in \go s_i}\lambda_te_t\big\|^2_{\go
s_i}\big)^{\frac12}\Del (\go s_i)_{i=1}^l \textrm { are incomparable
segments of  } \mathsf T\big\}
$$
and denote by $\ku T_2$ the completion of $\ku V$ under this norm.
The space $\ku T_2$ is what is called an $\ell_2$-Baire sum in
\cite{ad} and will play a universality role in the following. We
notice first that $(e_t)_{t\in\mathsf T}$ is a suppression
unconditional basis for $\ku T_2$, i.e., the projection onto any subsequence has norm $1$, and we need therefore not concern
ourselves with any particular enumeration of it in ordertype $\om$.
For any branch $\sigma=(\alpha,\beta)\in [\mathsf T]$, we denote by
$X_\sigma$ the closed subspace of $\ku T_2$ generated by the vectors
$(e_t)_{t\subseteq \sigma}$. Then $(e_t)_{t\subseteq
\sigma}=(e_{\sigma|n})_{n<\om}$ is a suppression unconditional
Schauder basis for $X_\sigma$, which by inequality (\ref{eq with p})
is $2$-equivalent with the standard unit vector basis of
$\ell_{p_\alpha}$.

We should note the following thing about the segment norm. Assume
$$
\go s=\{(u_0,s_0)\subsetneq (u_1,s_1)\subsetneq \ldots\subsetneq
(u_n,s_n)\}
$$
and
$$
\go s'=\{(u_0,s'_0)\subsetneq (u_1,s'_1)\subsetneq \ldots\subsetneq
(u_n,s'_n)\}
$$ are two segments of $\mathsf T$ whose first coordinates coincide,
then
\begin{equation}
\|\lambda_0e_{(u_0,s_0)}+\ldots+\lambda_ne_{(u_n,s_n)}\|_\go s
=\|\lambda_0e_{(u_0,s'_0)}+\ldots+\lambda_ne_{(u_n,s'_n)}\|_{\go
s'}.
\end{equation}
This will allow us to prove the following lemma.
\begin{lemme}\label{pres of norm}
Let $S,T$ be subtrees of  $\mathsf T$ and $\phi:S\til T$ an
isomorphism of trees preserving the first coordinates, i.e., for all
$(u,s)\in S$ there is some $s'$ such that $\phi(u,s)=(u,s')$. Then
the map
$$
M_\phi:e_{(u,s)}\mapsto e_{\phi(u,s)}
$$
extends to a surjective linear isometry from the space
$Z_S=[e_t]_{t\in S}\subseteq \ku T_2$ onto $Z_T=[e_t]_{t\in
T}\subseteq \ku T_2$.
\end{lemme}

\begin{proof}
By symmetry it suffices to prove that for any finite linear
combination $x$ of $(e_t)_{t\in S}$ we have $\tno x\tno \leq
\tno M_\phi(x)\tno $.

So fix incomparable segments $(\go s_i)_{i=1}^l$ of $\mathsf T$ and
consider the estimation
$$
\big (\sum_{i=1}^l\big\|\sum_{t\in \go s_i}\lambda_te_t\big\|^2_{\go
s_i}\big)^{\frac12}\leq \tno x\tno .
$$
As the support of $x$ is completely contained in $S$, we can, by
projecting onto suitable initial segments, suppose, without changing
the lower estimate of $\tno x\tno $, that each $\go s_i$ is completely
contained within $S$. But then
\begin{displaymath}\begin{split}
\big (\sum_{i=1}^l\big\|\sum_{t\in \go s_i}\lambda_te_t\big\|^2_{\go
s_i}\big)^{\frac12} = \big (\sum_{i=1}^l\big\|\sum_{t\in \go
s_i}\lambda_te_{\phi(t)}\big\|^2_{\phi[\go
s_i]}\big)^{\frac12}\leq\tno M_\phi(x)\tno .
\end{split}\end{displaymath}
So by taking suprema we see that $\tno x\tno \leq\tno M_\phi(x)\tno $.
\end{proof}

\section{Permutative equivalence}
We are now in a position to show that the relation of permutative
equivalence between (suppression unconditional) basic sequences is a
complete analytic equivalence relation.

We let $\ku {UBS}$ denote the standard Borel space of unconditional
basic sequences, i.e., $\ku {UBS}$ can be chosen to be the set of
subsequences of the universal unconditional basic sequence $(u_n)$
of Pe\l zcy\'nski (see \cite{p}). Denote by $(x_i)\approx_p(y_i)$
the fact that the two bases $(x_i)$ and $(y_i)$ in $\ku {UBS}$ are
permutatively equivalent, i.e., for some permutation $f$ of $\N$,
$(x_i)\approx (y_{f(i)})$. We recall that, as was first noticed by
Mityagin \cite{m}, unconditional basic sequences satisfy the
Schr\"oder--Bernstein principle, i.e., if $(x_i)$ and $(y_i)$ are
normalised unconditional basic sequences and $f,g:\N\til \N$
injections such that $(x_i)\approx(y_{f(i)})$ and
$(y_i)\approx(x_{g(i)})$, then $(x_i)$ and $(y_i)$ are permutatively
equivalent. This is easily seen to follow from the proof of the
Schr\"oder-Bernstein Theorem.

\begin{thm}
The relation of permutative equivalence, $\approx_p$, between
unconditional basic sequences is a complete analytic equivalence
relation.
\end{thm}

\begin{proof}
We shall reduce the pair $(\equiv_{\fed
\Sigma_1^1},\neq_{\fed\Sigma_1^1})$ between pruned normal trees on
$2\times \om$ to $\approx_p$. The reduction $\phi$ is the obvious
one given by
$$
\phi: S\mapsto (e_t)_{t\in S},
$$
where $(e_t)_{t\in S}$ is enumerated in ordertype $\om$ in some
canonical way. Here $(e_t)_{t\in\mathsf T}$ is the canonical basis
for the space $\ku T_2$. Since $(e_t)_{t\in S}$ is suppression
unconditional, it remains a basic sequence any way we enumerate it.

Suppose first that $S$ and $T$ are pruned normal trees on $2\times
\om$ such that $S\equiv_{\fed\Sigma_1^1}T$ as witnessed by some
$\alpha\in \om^\om$. Then, by Lemma \ref{pres of norm},
$e_{(u,s)}\mapsto e_{(u,s+\alpha|_{|s|})}$ restricts to an isometric
embedding of $[e_t]_{t\in S}$ into $[e_t]_{t\in T}$ and also to an
isometric embedding of $[e_t]_{t\in T}$ into $[e_t]_{t\in S}$. In
particular, the unconditional bases $(e_t)_{t\in S}$ and
$(e_t)_{t\in T}$ are equivalent to subsequences of each other and
hence permutatively equivalent.

On the other hand, if $S\neq _{\fed\Sigma_1^1}T$, we can find some
$\alpha\in 2^\om$ such that $\alpha\in {\rm proj}[S]\setminus{\rm
proj}[T]$. Take some $\beta\in \om^\om$ such that $(\alpha,\beta)\in
[S]$, and notice then that $(e_{(\alpha|n,\beta|n)})_n$ is
equivalent with the unit vector basis in $\ell_{p_\alpha}$. We claim
that there is no subsequence of $(e_t)_{t\in T}$ equivalent with
$\ell_{p_\alpha}$. To see this, notice that if $(e_t)_{t\in A}$ was
any subsequence of $(e_t)_{t\in T}$, then by Ramsey's theorem we could
find some infinite subset $B\subseteq A$ such that either
$B\subseteq \{(\gamma|n,\delta|n)\del n\in \N\}$ for some
$(\gamma,\delta)\in [T]$ or $B$ is an antichain in $T$.

In the first case, $(e_t)_{t\in B}$ is equivalent with a subsequence
of the unit vector basis of $\ell_{p_\gamma}$ and hence, as
$p_\gamma\neq p_\alpha$, is not equivalent with $\ell_{p_\alpha}$,
and in the latter case, by construction of $\ku T_2$, $(e_t)_{t\in
B}$ is equivalent with $\ell_2$, which is not equivalent with
$\ell_{p_\alpha}$ either. Thus $S\neq_{\fed\Sigma_1^1}T\saa
(e_t)_{t\in S}\not\approx_p(e_t)_{t\in T}$. This finishes the proof
of the reduction.
\end{proof}

We easily see from the above construction that we also reduce the
pair
$$
(\leq_{\fed\Sigma_1^1},\not\subseteq_{\fed\Sigma_1^1})
$$
to the
relation of being permutatively equivalent with a subsequence between
unconditional basic sequences. And thus
\begin{thm}
The relation between unconditional basic sequences of being
permutatively equivalent with a subsequence is a complete analytic quasiorder.
\end{thm}

The are several related results concerning equivalence of basic
sequences. For example, Ferenczi and Rosendal show in \cite{fr} that
a basic sequence is either subsymmetric, i.e., equivalent with all
its subsequences, or the relation $E_0$ Borel reduces to equivalence
between its subsequences. Also Rosendal \cite{r:cof} shows that the
relation of equivalence between basic sequences is Borel bireducible
with a complete $K_\sigma$ equivalence relation. And finally,
Ferenczi \cite{f} proves that if $(e_i)$ is an unconditional basic
sequence, then either $E_0$ Borel reduces to the relation of
permutative equivalence between the normalised blockbases of
$(e_i)$, or, for some $\ell_p$ or $c_0$, any normalised blockbasis
has a subsequence equivalent with this $\ell_p$ or $c_0$.

\section{Uniform homeomorphism of complete separable metric spaces}
We now intend to show that the relation of uniform homeomorphism
between complete separable metric spaces is a complete analytic
equivalence relation. This will be done by reducing the relation of
permutative equivalence between unconditional basic sequences to it.
Let us first remark that uniform homeomorphism is indeed analytic.
To see this, notice that if $(X,d_X)$ and $(Y,d_Y)$ are two complete
separable metric spaces, then they are uniformly homeomorphic if and
only if they have countable dense subsets $D_X$ and $D_Y$ that are
uniformly homeomorphic. For any uniform homeomorphism between $D_X$
and $D_Y$ will preserve Cauchy sequences in both directions and
hence extend to a uniform homeomorphism between $X$ and $Y$. But the
relation of uniform homeomorphism between countable metric spaces is
easily seen to be analytic, whence this extends to all complete
separable metric spaces.

This argument clearly does not extend to the relation of
homeomorphism between complete separable metric spaces (or more
naturally to the class of Polish topological spaces). A priori this
relation is not analytic but only $\fed \Sigma^1_2$, but, as we
shall see, it is $\fed\Sigma^1_1$-hard as an equivalence relation.
It is natural to ask the following.
\begin{quest}
Is the relation of homeomorphism between Polish spaces, i.e., closed
subspaces of $\R^\N$, a complete $\fed \Sigma^1_2$ equivalence
relation?
\end{quest}

In the following, we fix a normalised bimonotone unconditional basic
sequence $(e_n)_{n \in N}$ in a Banach space, where $N$ is an
unordered infinite countable set, and we let $X$ be the closed
subspace generated by the sequence $(e_n)_{n \in N}$, and denote by
$d_X$ the metric on $X$.

A {\em type} is a non-empty finite non-decreasing sequence of
strictly positive rational numbers. If
$t=(\lambda_1,\ldots,\lambda_n)$ is a type, we say that a
 vector $x \in
X$ {\em has type $t$} when $x$ can be written $x=\sum_{i=1}^m
\lambda_i e_{\sigma(i)}$, for some injection $\sigma$ of
$\{1,\ldots,m\}$ into $N$. Since $(e_n)$ is a basis, it is clear
that each vector of $X$ has at most one type. We enumerate the set
of types as $(t_n)_{n \in \N}$, with $t_1=(1)$, and we let $T_n$ be
the set of vectors of $X$ of type $t_n$; in particular $T_1$ is the
set of the unit vectors of the basis $(e_n)_{n \in N}$.

\begin{lemme}\label{separation} For any $n \in \N$, there exists $\delta_n>0$ such that the
set $T_n$ is $\delta_n$-separated.
\end{lemme}

\begin{proof} Write $t_n=(\lambda_1,\ldots,\lambda_m)$, let $\lambda_0=0$,
and let
$$\delta_n=\min(\{|\lambda_i-\lambda_j|\del 0 \leq i,j \leq m\} \cap \R_+^*).$$
If $x,y \in T_n$ are distinct, then there exists $k \in N$ such that
$p_k(x) \neq p_k(y)$, where $p_k$ denotes the (norm $1$) projection
onto $[e_k]$. Since $p_k(x)$ (respectively $p_k(y)$) is equal to
$\lambda e_k$ for some $\lambda \in \{\lambda_i\del 0 \leq i \leq
m\}$, it follows that $|p_k(x)-p_k(y)| \geq \delta_n$, and therefore
that $\norm{x-y} \geq |p_k(x-y)| \geq \delta_n.$
\end{proof}

We now describe how to build a Polish space $P(X)$ by implanting on
$X$ various elementary metric spaces in order to rigidify its
topological structure.

For any $n \in \N$, let $H_n$ be a fixed metric space with a special
point $0_n$, which is the union of $n$ isometric copies of $[0,1]$,
each of which has $0_n$ as endpoint, and which intersect only in
$0_n$. Now for any $n \in \N$ and $x \in X$ of type $t_n$, we let
$H(x)$ be an isometric copy of $H_n$ in which we denote the special
point corresponding to $0_n$ by $x$, and we denote by $d_x$ the
metric on $H(x)$. We also denote $H^0(x)=H(x) \setminus \{x\}$.

We then let $P(X)$ be the amalgamation of $X$ with all $H(x)$, for
$n \in \N$ and $x \in T_n$, each $H(x)$ being amalgamated with
$P(X)$ in $x$. This means that
$$
P(X)=X \sqcup (\bigsqcup_{n \in \N}\bigsqcup_{x \in T_n} H^0(x)),
$$
where the metric $d$ on $P(X)$ is defined as follows, for $y,z$ in
$P(X)$:

- if $y$ and $z$ both belong to $X$, then $d(y,z)=d_X(y,z)$,

- if $y$ and $z$ both belong to some $H^0(x)$, then
$d(y,z)=d_x(y,z)$,

- if $y \in X$ and $z$ belongs to some $H^0(x)$, then
$d(y,z)=d_X(y,x)+d_x(x,z)$,

- if $y$ belongs to some $H^0(x)$ and $z$ belongs to some $H^0(x')$,
with $x \neq x'$, then $d(y,z)=d_x(y,x)+d_X(x,x')+d_{x'}(x',z)$.

The set $R=\bigcup_{n \in \N}T_n \subset P(X)$ is called the set of
{\em roots} in $P(X)$, and the set $H=\bigcup_{x \in R} H^0(x)
\subset P(X)$ is called the {\em hair} in $P(X)$.  We have the
following fact:

\begin{lemme}
The space $P(X)$ is separable and complete metric.
\end{lemme}

\begin{proof}
$P(X)$ is obviously separable. If $(y_k)_{k \in \N}$ is a Cauchy
sequence in $P(X)$, and $y_k$ belongs to $X$ for all $k \in \N$,
then $(y_k)_k$ converges in $X$; therefore we may assume that $y_k$
belongs to the hair for all $k \in \N$. If there is a fixed $x \in
X$ such that $y_k$ belongs to  $H^0(x)$ for all $k \in \N$, then
$(y_k)_k$ converges in $H^0(x) \cup \{x\}$ by completeness of
$H(x)$; therefore we may assume that there is a sequence $(x_k)_{k
\in \N}$ of pairwise distinct points of $X$ such that $y_k \in
H^0(x_k)$ for all $k \in \N$. Then for all $j,k$ in $\N$,
$$d(y_j,y_k)=d(y_j,x_j)+d(x_j,x_k)+d(x_k,y_k),$$ and we deduce that $(x_k)_k$ is a Cauchy
sequence in $X$ and that $(d(x_k,y_k))_k$ converges to $0$;
therefore the sequence $(y_k)_k$ converges to some $x$ in $X$.
\end{proof}

\begin{prop}
Let $(e_n)_{n \in N}$ and $(e'_n)_{n \in N'}$ be normalised
bimonotone unconditional basic sequences  and let $X=[e_n, n \in
N]$, $X'=[e'_n, n \in N']$. Then any homeomorphism between $P(X)$
and $P(X')$ takes $X$ onto $X'$, the hair in $P(X)$
 onto the hair in $P(X')$, the set of roots in $P(X)$ onto the set of roots in $P(X')$,
and for each $n \in \N$, the set of type $t_n$ vectors of $X$ onto
the set of type $t_n$ vectors of $X'$.
\end{prop}

\begin{proof}
Indeed $H$ is the set of points in $P(X)$ which admit an open
neighborhood homeomorphic to $[0,1[$ or $]0,1[$, and $X=P(X)
\setminus H$. Also, define an {\em implant} as a maximal subset of
$H$ homeomorphic to $[0,1[$, and given $x \in X$, say that an
implant $h$ is {\em attached to} $x$ if $x$ is adherent to $h$. It
is then clear that a point $x$ in $P(X)$ is a root if and only if
some implant is attached to it, and that $x$ is a point of type
$t_n$, for $n \in \N$, if and only if exactly $n$ implants are
attached to $x$.
\end{proof}

\begin{prop}\label{uniquehomeo}
Let $(e_n)_{n \in N}$ and $(e'_n)_{n \in N'}$ be normalised
bimonotone unconditional basic sequences  and let $X=[e_n, n \in
N]$, $X'=[e'_n, n \in N']$. Let $T$ be a homeomorphism between
$P(X)$ and $P(X')$. Then there exists a bijection $\sigma$ between
$N$ and $N'$ such that for any finite subset $I$ of $N$, any
sequence $(\lambda_i)_{i \in I}$ of non-negative real numbers,
$$T(\sum_{n \in I}\lambda_n e_n)= \sum_{n \in I}\lambda_n
e^{\prime}_{\sigma(n)}.$$
\end{prop}

\begin{proof} Since $T$ maps type $(1)$ points of $P(X)$ onto type $(1)$
points of $P(X')$, there exists a bijection $\sigma$ between $N$ and
$N'$ such that $T(e_n)=e'_{\sigma(n)}$ for all $n \in N$. By
continuity of $T$, it is then enough to prove by induction on $|I|$
that for any finite subset $I$ of $N$, and for any sequence
$(\lambda_n)_{n \in I}$ of pairwise distinct positive rationals,
$$
T(\sum_{n \in I}\lambda_n e_n)=\sum_{n \in I}\lambda_n
e'_{\sigma(n)}.
$$

For any $n \in \N$, any $\lambda \in \Q^{+*}$, $T(\lambda e_n)$ has
type $(\lambda)$, therefore $T(\lambda e_n)=\lambda
e'_{k_n(\lambda)}$ for some $k_n(\lambda)$ in $N'$. By continuity,
$e'_{k_n(\lambda)}=T(\lambda e_n)/\lambda$ is constant on $\Q^{+*}$
and equal to $e'_{k_n(1)}=T(e_n)=e'_{\sigma(n)}$. Therefore
$$T(\lambda e_n)=\lambda e'_{\sigma(n)}, \forall \lambda \in \Q^{+*}.$$

Let now $I \subset N$ be finite, with $|I| \geq 2$. Let $\Delta$ be
the open subset of $(\Q^{+*})^I$ defined by
$$\Delta=\{(\lambda_n)_{n \in I}: \forall n \neq p, \lambda_n \neq \lambda_p\}.$$
For $\lambda=(\lambda_n)_{n \in I} \in \Delta$, let
$x(\lambda)=\sum_{n \in I}\lambda_n e_n$. Then $T(x(\lambda))$ has
the same type as $x(\lambda)$ and therefore may be written
(uniquely) in the form
$$
T(x(\lambda))=\sum_{n \in I}\lambda_n e'_{k_n(\lambda)},
$$
with $k_n(\lambda) \in N'$ for each $n \in I$.

We prove that $k_n(\lambda)$ is locally constant on $\Delta$, for
all $n \in I$. Fix indeed $\lambda \in \Delta$, then for any
$\mu=(\mu_n)_{n \in I}$ in a neighborhood  $V$ of $\lambda$ in
$\Delta$, we have
\begin{displaymath}\begin{split}
&\norm{\sum_{n \in I} \lambda_n e'_{k_n(\lambda)}- \sum_{n \in I}
\lambda_n e'_{k_n(\mu)}}\\
 \leq &\norm{\sum_{n \in I} \lambda_n
e'_{k_n(\lambda)}- \sum_{n \in I} \mu_n e'_{k_n(\mu)}}+\norm{\sum_{n
\in
I}(\lambda_n-\mu_n)e'_{k_n(\mu)}}\\
\leq& \norm{T(x(\lambda))-T(x(\mu))}+\sum_{n \in
I}|\lambda_n-\mu_n|.
\end{split}\end{displaymath}
Therefore if $V$ is small enough,
$$\norm{\sum_{n \in I} \lambda_n e'_{k_n(\lambda)}-
\sum_{n \in I} \lambda_n e'_{k_n(\mu)}}<\delta, \forall \mu \in V,$$
where $\delta$ is such that the set of points of $X$ of same type as
$x$ is $\delta$-separated (Lemma \ref{separation}). Therefore,
$$\sum_{n \in I}\lambda_n e'_{k_n(\lambda)}=
\sum_{n \in I}\lambda_n e'_{k_n(\mu)}, \forall \mu \in V,$$ and
since the $\lambda_n$'s for $n$ in $I$ are pairwise distinct,
$$k_n(\lambda)=k_n(\mu),\ \forall n \in I, \forall \mu \in V.$$

We deduce from this fact that for all $n$ in $I$, $k_n$ is constant
on each connected component of $\Delta$.

Fix now $\lambda \in \Delta$, let $C(\lambda)$ be the connected
component of $\Delta$ containing $\lambda$, and let $(k_n)_{n \in I}
\in (N')^I$ be such that
$$T(x(\lambda))=\sum_{n \in I} \lambda_n e'_{k_n}.$$
Let $n_0 \in I$ be such that $\lambda_{n_0}=\min_{n \in I}
\lambda_n$. For any $t \in ]0,\lambda_{n_0}[$, the element of
$\Delta$ associated to $\sum_{n \in I, n \neq n_0} \lambda_n e_n+t
e_{n_0}$ is in $C(\lambda)$, therefore
$$
T(\sum_{n \in I, n \neq n_0} \lambda_n e_n+t e_{n_0})
=\sum_{n \in I, n \neq n_0} \lambda_n e'_{k_n}+te'_{k_{n_0}}.
$$
When
$t$ converges to $0$, we obtain that
\begin{displaymath}\begin{split}
\sum_{n \in I, n \neq n_0} \lambda_n e'_{k_n}=&T(\sum_{n \in I, n
\neq n_0} \lambda_n e_n)\\
=&\sum_{n \in I, n \neq n_0} \lambda_n e'_{\sigma(n)},
\end{split}\end{displaymath}
by the induction hypothesis. Since the $\lambda_n$'s, $n \in I$, are
pairwise distinct, it follows that $\a n\in I\setminus \{n_0\}$,
$k_n=\sigma(n)$.

Let now $n_1 \in I$ be such that $\lambda_{n_1}=\min_{n \neq n_0}
\lambda_n$ and let $\lambda^{sym} \in \Delta$ be defined by
$\lambda^{sym}_n=\lambda_n,\ \forall n \notin \{n_0,n_1\}$,
$\lambda^{sym}_{n_0}=\lambda_{n_1}$, and
$\lambda^{sym}_{n_1}=\lambda_{n_0}$. There exists $(l_n)_{n \in I}
\in (N')^I$ such that for any $\mu=(\mu_n)_{n \in I}$ in
$C(\lambda^{sym})$,
$$
T(\sum_{n \in I} \mu_n e_n)=\sum_{n \in I} \mu_n e'_{l_n},
$$
and the same reasoning as above gives us that
$$
l_n=\sigma(n),\ \forall n \in I, n \neq n_1.
$$
Now
\begin{displaymath}\begin{split}
T(\sum_{n \neq n_0} \lambda_n e_n+\lambda_{n_1}e_{n_0})=& \lim_{t
\rightarrow (\lambda_{n_1})^{-}} T(\sum_{n \neq n_0} \lambda_n e_n
+t e_{n_0})\\
=&\sum_{n \neq n_0} \lambda_n e'_{k_n} + \lambda_{n_1} e'_{k_{n_0}},
\end{split}\end{displaymath}
since the element of $\Delta$ associated to $\sum_{n \neq
n_0}\lambda_n e_n +te_{n_0}$ is in $C(\lambda)$ for each $t$ in
$[\lambda_{n_0},\lambda_{n_1}[$. Also,
 $$T(\sum_{n \neq n_0} \lambda_n e_n+\lambda_{n_1}e_{n_0})=
\lim_{t \rightarrow (\lambda_{n_1})^{-}} T(\sum_{n \neq n_1}
\lambda_n^{sym} e_n +t e_{n_1})$$
$$=\sum_{n \neq n_1} \lambda_n^{sym} e'_{l_n}+\lambda_{n_1}e'_{l_{n_1}}=
\sum_{n \neq n_0} \lambda_n e'_{l_n}+\lambda_{n_1} e'_{l_{n_0}},$$
since the element of $\Delta$ associated to $\sum_{n \neq
n_1}\lambda_n^{sym} e_n +te_{n_1}$ is in $C(\lambda^{sym})$ for each
$t \in [\lambda_{n_0},\lambda_{n_1}[$. Therefore
$$\sum_{n \neq n_0} \lambda_n e'_{k_n} + \lambda_{n_1} e'_{k_{n_0}}
=\sum_{n \neq n_0} \lambda_n e'_{l_n} + \lambda_{n_1}
e'_{l_{n_0}},$$ from which it follows that
$$\{k_{n_0},k_{n_1}\}=\{l_{n_0},l_{n_1}\}.$$
Since $k_{n_1}=\sigma(n_1)$ and $l_{n_0}=\sigma(n_0)$, we deduce
that $k_{n_0}=\sigma(n_0)$.

We have finally proved that $k_n=\sigma(n)$, $\forall n \in I$, and
therefore
$$T(\sum_{n \in I} \lambda_n e_n)=\sum_{n \in I} \lambda_n e'_{\sigma(n)}, \ \forall
(\lambda_n)_{n \in I} \in \Delta.$$
\end{proof}

\begin{thm} The pair (Lipschitz isomorphism, non-homeomorphism) restricted to the class of complete separable metric spaces  is
hard for analytic equivalence relations, and thus, the relation of
uniform homeomorphism between complete separable metric spaces is a
complete analytic equivalence relation.
\end{thm}

\begin{proof} We show that the relation of permutative equivalence between
subsequences of the universal unconditional sequence of Pe\l
czy\'nski $(u_n)_{n \in \N}$ is reducible to the pair (Lipschitz isomorphism, non-homeomorphism) restricted to the class of complete separable metric spaces. Up to equivalent renorming, we may assume that $(u_n)_{n \in
\N}$ is bimonotone.

For $N$ an infinite subset of $\N$, we define $$\alpha(N)=P([u_n, n
\in N]).$$ For any $N$, $\alpha(N)$ is canonically isometric to a
closed subset of $P([u_n, n \in \N])$. In this setting, the map
$\alpha$ is clearly Borel. Furthermore, whenever $(u_n)_{n \in N}$
and $(u_n)_{n \in N'}$ are permutatively equivalent, there is a
natural Lipschitz isomorphism between $\alpha(N)$
and $\alpha(N')$.

Conversely, if $T$ is a homeomorphism between $\alpha(N)$ and
$\alpha(N')$, then by Proposition \ref{uniquehomeo}, there exists a
bijection $\sigma$ between $N$ and $N'$ such that for any finite
subset $I$ of $N$, for any sequence $(\lambda_n)_{n \in I}$ of
non-negative reals,
$$
T(\sum_{n \in I}\lambda_n u_n)=\sum_{n \in I}\lambda_n u_{\sigma(n)}.
$$
It follows that $(u_n)_{n \in N} \approx_p (u_n)_{n \in N'}$.
Indeed, let $(\lambda_n)_n \in \R^{N}$ be such that $\sum_{n \in N}
\lambda_n u_n$ converges, then $\sum_{n \in N} |\lambda_n| u_n$
converges by unconditionality, and therefore
$$
\sum_{n \in N}
|\lambda_n|u_{\sigma(n)}=T(\sum_{n \in N} |\lambda_n| u_n)
$$
converges, so $\sum_{n \in \N} \lambda_n u_{\sigma(n)}$ converges,
again by unconditionality. Conversely,
$$
\sum_{n \in N} \lambda_n
u_n
$$
converges whenever
$$
\sum_{n \in N} \lambda_n u_{\sigma(n)}
$$
converges. We deduce that $(u_n)_{n \in N}$ is equivalent to
$(u_{\sigma(n)})_{n \in N}$.
\end{proof}

Unfortunately, we have not been able to replace our spaces $P(X)$ by
Banach spaces and thus the following problem remains open.
\begin{prob}
What is the complexity with respect to $\leq_B$ of the relation of
uniform homeomorphism between separable Banach spaces? In
particular, is it analytic complete?\end{prob}

The corresponding quasiorder of uniform homeomorphic embeddability
has also been studied in the form of homeomorphic embeddability
between compact metric spaces. A series of results by
Louveau--Rosendal \cite{lr}, Marcone--Rosendal \cite{mr}, and
culminating in Camerlo \cite{c}, show that the relation of
continuous embeddability between dendrites all of whose branching
points have order three is a complete analytic quasiorder.

\section{Isomorphism of separable Banach spaces}
In order to prove our main result that isomorphism of separable
Banach spaces is complete, the simple construction of the spaces
$Z_S$ does not seem to suffice. For example, it is known (see, e.g.,
\cite{ad}) that the space $\ku T_2$ contains a copy of $c_0$ and
therefore the control over the subspaces present is presumably not
good enough. Instead, we shall use the Davis, Figiel, Johnson, and
Pe\l czy\'nski \cite{dfjp} interpolation method and the results
proved in \cite{ad} to avoid certain subspaces.

Suppose $S$ is a pruned subtree of $\mathsf T$. We denote by $Z_S$
the closed subspace of $\ku T_2$ spanned by $(e_t)_{t\in S}$. The
latter is still a suppression unconditional basis for $Z_S$.

In order to obtain a better control of the subspaces present, we
shall now replace $Z_S$ with an interpolate that eliminates some
vectors whose support is too much in between several different
branches.
\begin{defi}
Let $W_S$ be the convex hull in $Z_S$ of the set $\bigcup_{\sigma\in
[S]}B_{X_\sigma}$ and let for each $n\geq 0$, $C_S^n$ be the convex
set $2^nW_S+2^{-n}B_{Z_S}$. As $W$ is a bounded set, we can for each
$n$ define an equivalent norm on $Z_S$ by the taking the gauge of
$C_S^n$
$$
\|x\|_S^n:=\inf (\lambda\del \frac x\lambda\in C_S^n).
$$
\end{defi}
Our first lemma shows that the $n$'th norm of a vector does not
depend on the ambient space.
\begin{lemme}
Let $S$ and $T$ be two pruned trees and $x\in \ku V$ a finitely
supported vector belonging to both $Z_S$ and $Z_T$. Then for every
$n$,
$$
\|x\|_S^n=\|x\|_T^n.
$$
\end{lemme}

\begin{proof}
Obviously, by symmetry, it is enough to prove that
$\|x\|_S^n\geq\|x\|_T^n$. So suppose $\lambda>0$ is such that $\frac
x\lambda\in C^n_S=2^nW_S+2^{-n}B_{Z_S}$. Then we can find a finite
number of branches $\sigma_1,\ldots,\sigma_m\in [S]$, vectors
$y_i\in B_{X_{\sigma_i}}$, scalars $r_1,\ldots,r_m>0$, and $z\in
B_{Z_S}$ such that $\sum_ir_i=1$ and
$$
\frac x\lambda=2^n(r_1y_1+\ldots+r_my_m)+2^{-n}z.
$$
Since $x\in Z_T$ and has finite support, we can choose a finite
number of branches $\chi_1,\ldots,\chi_k\in [T]$ such that ${\rm
support}(x)\subseteq \chi_1\cup\ldots\cup\chi_k$. Let $R$ be the
pruned tree whose branches are $\chi_1,\ldots,\chi_k$ and let $P_R$
be the canonical projection of $\ku T_2$ onto $[e_t]_{t\in R}$. As
$(e_t)_{t\in \mathsf T}$ is suppression unconditional, $\|P_R\|=1$.
Thus, $P_R(z)\in B_{Z_T}$ and, as the $y_i$ belong to subspaces
spanned by branches, for each $i=1,\ldots,m$ there is some $1\leq
i'\leq k$, such that $P_R(y_i)\in B_{X_{\chi_{i'}}}$. Therefore,
\begin{displaymath}\begin{split}
\frac x\lambda=&P_R(\frac
x\lambda)\\
=&P_R(2^n(r_1y_1+\ldots+r_my_m)+2^{-n}z)\\
=&2^n(r_1P_R(y_1)+\ldots+r_mP_R(y_m))+2^{-n}P_R(z)\\
\in& 2^nW_T+2^{-n}B_{Z_T}=C^n_T.
\end{split}\end{displaymath}
And hence,
$$
\|x\|_S^n=\inf (\lambda\del \frac x\lambda\in C_S^n)\geq\inf
(\lambda\del \frac x\lambda\in C_T^n)=\|x\|_T^n.
$$
\end{proof}

\begin{lemme}
Suppose $\phi:S\til T$ is an isomorphism of pruned subtrees of
$\mathsf T$ satisfying $\phi(u,s)=(u,s')$ , i.e., $\phi$ preserves
the first coordinate of every element of $S$. Then for every $n$,
the mapping
$$
M_\phi:e_{(u,s)}\mapsto e_{\phi(u,s)}
$$
extends (uniquely) to a surjective linear isometry from
$(Z_S,\|\cdot\|_S^n)$  onto $(Z_T,\|\cdot\|_T^n)$.
\end{lemme}

\begin{proof}
By symmetry it is again enough to show that for any finitely
supported vector $x\in Z_S$ we have
$\|x\|_S^n\geq\|M_\phi(x)\|_T^n$. So suppose $\lambda>0$ is such
that $\frac x\lambda\in C^n_S=2^nW_S+2^{-n}B_{Z_S}$, and find a
finite number of branches $\sigma_1,\ldots,\sigma_m\in [S]$, vectors
$y_i\in B_{X_{\sigma_i}}$, scalars $r_1,\ldots,r_m>0$, and $z\in
B_{Z_S}$ such that $\sum_ir_i=1$ and
$$
\frac x\lambda=2^n(r_1y_1+\ldots+r_my_m)+2^{-n}z.
$$
But then by Lemma \ref{pres of norm}, $M_\phi(y_i)\in
B_{X_{\phi[\sigma_i]}}$ for each $i$, while, as $M_\phi$ is an
isometry from $(Z_S,\tno \cdot\tno )$ to $(Z_T,\tno \cdot\tno )$, also
$M_\phi(z)\in B_{Z_T}$. Hence, $\frac{M_\phi(x)}\lambda\in C_T^n$
and $\|x\|_S^n\geq\|M_\phi(x)\|_T^n$.
\end{proof}

Combining the two preceding lemmas we have
\begin{lemme}
Suppose $\phi:S\til T$ is an embedding of pruned subtrees of
$\mathsf T$ satisfying $\phi(u,s)=(u,s')$ , i.e., $\phi$ preserves
the first coordinate of every element of $S$. Then for every $n$,
the mapping
$$
M_\phi:e_{(u,s)}\mapsto e_{\phi(u,s)}
$$
extends to a linear isometry from $(Z_S,\|\cdot\|_S^n)$ into
$(Z_T,\|\cdot\|_T^n)$.
\end{lemme}

\begin{defi}
Let $S$ be a pruned subtree of $\mathsf T$ and let
$\ell_2(Z_S,\|\cdot\|_S^n)$ be the $\ell_2$-sum of the sequence of
spaces $(Z_S,\|\cdot\|_S^n)_n$. We denote by $\Delta(Z_S,2)$ the
closed subspace of $\ell_2(Z_S,\|\cdot\|_S^n)$ consisting of all the
vectors of $\ell_2(Z_S,\|\cdot\|_S^n)$ on the form $(x,x,x,\ldots)$.
\end{defi}

We should note that, as for each $t\in S$, $e_t\in W_S$, we have
$\|e_t\|_S^n\leq 2^{-n}$, and hence $(e_t,e_t,e_t,\ldots)\in
\Delta(Z_S,2)$. As $(e_t)_{t\in S}$ remains a suppression
unconditional basis for $(Z_S,\|\cdot\|_S^n)$ for each $n$, one also
sees that it is a suppression unconditional basis for
$\Delta(Z_S,2)$. We shall denote by $\|\cdot\|_S$ the norm on
$(e_t)_{t\in S}$ giving the space $\Delta(Z_S,2)$.

\begin{prop}\label{embedding}
Suppose $\phi:S\til T$ is an embedding of pruned subtrees of
$\mathsf T$ satisfying $\phi(u,s)=(u,s')$ , i.e., $\phi$ preserves
the first coordinate of every element of $S$. Then
$$
M_\phi:e_{(u,s)}\mapsto e_{\phi(u,s)}
$$
extends to a linear isometry from $\Delta(Z_S,2)$ into
$\Delta(Z_T,2)$. Moreover, $M_\phi(\Delta(Z_S,2))$ is
$1$-complemented in $\Delta(Z_T,2)$.
\end{prop}

We now need the following fundamental result of Argyros and Dodos on
the structure of the spaces $\Delta(Z_S,2)$. We shall formulate
their result only for the special case of the spaces that we
construct here, which are particular examples of the more general
construction in \cite{ad}, and only mention the aspects we need.
\begin{thm}[S. Argyros and P. Dodos \cite{ad}, Theorem 71, Theorem 74]\label{argyros}
Let $S$ be a pruned subtree of the complete tree on $2\times \om$,
$\mathsf T$. For each $\sigma\in [S]$, denote by $\ku X_\sigma$ the
closed subspace of $\Delta(Z_S,2)$ spanned by the sequence
$(e_t)_{t\subseteq \sigma}$ and by $P_\sigma$ the (norm $1$)
projection of $\Delta(Z_S,2)$ onto $\ku X_\sigma$.
\begin{enumerate}
  \item For each $\sigma\in [S]$, $\ku X_\sigma$ is isomorphic to
  $X_\sigma\subseteq Z_S$.
  \item If $Y\subseteq \Delta(Z_S,2)$ is an infinite-dimensional
  closed subspace such that for all closed infinite-dimensional
  subspaces $Z\subseteq Y$ and $\sigma\in [S]$ the projection
  $P_\sigma:Z\til \ku X_\sigma$ is not an isomorphic embedding (in this case we
  say that $Y$ is $Z_S$-singular), then $Y$ contains $\ell_2$.
 \end{enumerate}
Moreover, $\Delta(Z_S,2)$ is reflexive.
\end{thm}

We are now ready to prove the main result of this article.
\begin{thm}
The relation of isomorphism between separable Banach spaces is a
complete analytic equivalence relation.
\end{thm}

Our proof will at the same time also show the following two results
\begin{thm}\label{lipschitz isomorphism}
The relation of Lipschitz isomorphism between separable Banach
spaces is a complete analytic equivalence relation.
\end{thm}
and
\begin{thm}
The relations of embeddability, complemented embeddability and
Lipschitz embeddability between separable Banach spaces are complete
analytic quasiorders.
\end{thm}
For good order, we should mention that by Theorem \ref{lipschitz
isomorphism}, the first problem of \cite{r:cof} is answered. Theorem
\ref{lipschitz isomorphism} should also be contrasted with the
result in \cite{r:cof} stating that the relation of Lipschitz
isomorphism between compact metric spaces is Borel bireducible with
a complete $\fed K_\sigma$-equivalence relation. Thus Lipschitz
isomorphism between compact metric spaces has the same complexity as
equivalence between  Schauder bases, while between separable Banach
spaces it has the same complexity as permutative equivalence.

\begin{proof}
The map that will simultaneously take care of all the reductions is
the obvious one
$$
S\mapsto \Delta(Z_S,2)
$$
for all pruned normal subtrees $S$ of $\mathsf T$.

We thus only need to notice the properties of this map. First of
all, if $S$ and $T$ are two pruned normal trees such that
$S\leq_{\fed\Sigma^1_1}T$ as witnessed by some $\beta\in \om^\om$,
then we can define an embedding $\phi:S\til T$ by
$\phi(u,s)=(u,s+\beta|_{|s|})$.  By Proposition \ref{embedding}, the
map $M_\phi:\Delta(Z_S,2)\til \Delta(Z_T,2)$ is an isomorphic
embedding, which moreover is a permutative equivalence between $(e_t)_{t\in S}$
and a subsequence of $(e_t)_{t\in T}$. So, in particular, if
$S\equiv_{\fed\Sigma_1^1}T$, then $(e_t)_{t\in S}$ and $(e_t)_{t\in
T}$ are permutatively equivalent with subsequences of each other and hence, as
they are both unconditional, they are permutatively equivalent and thus
$\Delta(Z_S,2)$ and $\Delta(Z_T,2)$ are isomorphic.

On the other hand, if $S\not\subseteq_{\fed\Sigma_1^1}T$, we find an
$\alpha\in {\rm proj}[S]\setminus{\rm proj}[T]$ and a $\beta$ such
that $\sigma=(\alpha,\beta)\in [S]$. We thus notice that
$\ell_{p_\alpha}\iso X_\sigma\iso\ku X_\sigma$. We claim that
$\Delta(Z_T,2)$ contains no subspace isomorphic to
$\ell_{p_\alpha}$. For if $Y$ is any subspace of $\Delta(Z_T,2)$,
then either $Y$ is $Z_T$-singular, in which case, $Y$ contains a
copy of $\ell_2$ and hence is not isomorphic to $\ell_{p_\alpha}$,
or there is a subspace $Z\subseteq Y$ and a branch
$\rho=(\gamma,\delta)\in[T]$ such that $P_\rho:Z\til \ku X_\rho$ is
an isomorphic embedding. But then $Z$ is isomorphic to a subspace of
$\ku X_\rho\iso \ell_{p_\rho}$ and hence contains a subspace
isomorphic to $\ell_{p_\rho}$. As $\rho\neq\alpha$, $Y$ cannot be
isomorphic to $\ell_{p_\alpha}$, and thus, finally, $\Delta(Z_S,2)$
does not embed into $\Delta(Z_T,2)$. Since by  Theorem \ref{argyros} $\Delta(Z_S,2)$ is reflexive, it follows that $\Delta(Z_S,2)$ does not Lipschitz embed into $\Delta(Z_T,2)$ (see \cite{bl} chapter 7 for more on this).

This shows that
$(\leq_{\fed\Sigma_1^1},\not\subseteq_{\fed\Sigma_1^1})$ reduces to
the couple (complemented isomorphic embeddability, non-Lipschitz embeddability) between
separable Banach spaces and thus the relations of complemented embeddability, embeddability, and Lipschitz embeddability are analytic
complete. Similarly,
$(\equiv_{\fed\Sigma_1^1},\neq_{\fed\Sigma_1^1})$ reduces to the
relations of isomorphism and Lipschitz isomorphism and these are analytic complete too.
\end{proof}

\begin{cor}
The relations of topological embeddability and topological
isomorphism between Polish groups are complete analytic as
quasiorders and equivalence relations respectively.
\end{cor}

Before we prove this, let us first define the space of Polish groups,
$\go G$, as the Effros--Borel space of closed subgroups of ${\rm
Hom}([0,1]^\N)$. By a result of Uspenski\u\i{} \cite{usp}, this
group contains all other Polish groups as closed subgroups up to
topological isomorphism. Two Polish groups are said to be {\em
topologically isomorphic} if there is a continuous group isomorphism
between them. Such an isomorphism is automatically a homeomorphism
and thus an isomorphism of the corresponding uniform structures.
Similarly, one Polish group is {\em topologically embeddable} into
another if it is topologically isomorphic with a closed subgroup.

\begin{proof}
Notice that if $\phi:X\til Y$ is a topological isomorphism of two
Banach spaces considered as Polish groups, then, in particular,
$\phi$ is an isomorphism of $X$ and $Y$ as $\Q$-vector spaces (since
it preserves divisibility).  But any continuous $\Q$-vector space
isomorphism between two Banach spaces is also a linear isomorphism.
The same argument applies to embeddings. So group
isomorphism/embedding coincides with linear isomorphism/embedding.
\end{proof}

\

\

\begin{flushleft}
{\em Address of V. Ferenczi and A. Louveau:}\\
\'Equipe d'Analyse Fonctionnelle,\\
Universit\'e Pierre et Marie Curie - Paris 6,\\
Bo\^ite 186,\\
4, Place Jussieu,\\
75252 Paris Cedex 05, France.\\
\texttt{ferenczi@ccr.jussieu.fr, louveau@ccr.jussieu.fr}
\end{flushleft}

\

\begin{flushleft}
{\em Address of C. Rosendal:}\\
Department of mathematics,\\
University of Illinois at Urbana-Champaign,\\
273 Altgeld Hall, MC 382,\\
1409 W. Green Street,\\
Urbana, IL 61801, USA.\\
\texttt{rosendal@math.uiuc.edu}
\end{flushleft}

\end{document}